% Logic Eprints
%Submitted 1021 Wed Jul 05, 1995 by: gitik@math.tau.ac.il (gitik moti )
%logic/gitik/addcohen795.tex
%

%Copied from Diana 4-7-95
%Corrected 4-7-95
% Logic Eprints
%Submitted 0312 Wed Mar 01, 1995 by: gitik@math.tau.ac.il (gitik moti )
%logic/gitik-shelah/lns
%logic/gitik-shelah/lns/def1200.tex
%

%This file is called def1200.tex you will need it
%to process
%the other file(s) you receive from me.

%-------------------------------------------------------
\def\today{\ifcase\month\or January\or February\or
March\or April\or May\or June\or July\or August\or
September\or October\or November\or December\fi
\space\number\day, \number\year}

%\font\lcmss=lcmssb8

%\font\cmss=cmssq10 at 10pt

%\font \rrrm=cmbx10 at 12pt

\def\dspace{\lineskip=2pt\baselineskip=18pt
\lineskiplimit=0pt}

\font \bbrm=cmbx10 at 12pt

\def\bigtype{\bbrm}

\hsize=13.5cm
\magnification=1200
\def\ce{\centerline}

\def\hb{\hfill\break}

\def\title #1{\null\bigskip\ce{\bigtype #1}
\bigskip}
%for discretionary break in equation use

%Greek letters
\def\alp{\alpha}		
\def\bet{\beta}		
\def\gam{\gamma}		
\def\del{\delta}

\def\kap{\kappa}
\def\lam{\lambda}

\def\ome{\omega}		

%Caligraphic roman letters

\def\calP{{\cal P}}

%Bold roman letters

%bold Greek capital letters

%Capital roman double letters
%\def\CC{{\rlap {\raise 0.4ex \hbox{$\scriptscriptstyle
%|$}}
%\hskip -0.15em C}}
%\def\AA{{\mathchoice
%{I\hskip -3.7pt {\rm A}}
%{I\hskip -3.7pt {\rm A}}
%{I\hskip -3.1pt {\rm A}}
%{I\hskip -2.5pt {\rm A}}}}
%\def\BB{{I\!\!B}}
%\def\EE{{I\!\!E}}
%\def\FF{{I\!\!F}}
%\def\NN{{I\!\!N}}
%\def\PP{{I\hskip-2.5pt P}}
%\def\QQ{{\rlap {\raise 0.4ex
%\hbox{$\scriptscriptlsubspacesstyle |$}}
%\hskip -0.1em Q}}
%\def\RR{{I\!\!R}}
%\def\ZZ{{Z\!\!\! Z}}
    
%Special fonts
\font\tenboldgreek=cmmib10
 \font\sevenboldgreek=cmmib10 at 7pt
\font\fiveboldgreek=cmmib10 at 7pt
\newfam\bgfam
\textfont\bgfam=\tenboldgreek
\scriptfont\bgfam=\sevenboldgreek
\scriptscriptfont\bgfam=\fiveboldgreek

%\def\bg{\fam6}
%\mathchardef\alpha="700B
%\def\bfalp{{\fam=\bgfam\balp}}
\mathchardef\ggarrow="7010

\font\tengerman=eufm10 \font\sevengerman=eufm7
\font\fivegerman=eufm5
\font\tendouble=msym10 \font\sevendouble=msym7
\font\fivedouble=msym5

\textfont4=\tengerman \scriptfont4=\sevengerman
\scriptscriptfont4=\fivegerman
\newfam\dbfam
\textfont\dbfam=\tendouble \scriptfont\dbfam=
\sevendouble
\scriptscriptfont\dbfam=\fivedouble

\mathchardef\ng="702D
\mathchardef\dbA="7041
\mathchardef\sm="7072
\mathchardef\nvdash="7030
\mathchardef\nldash="7031
\mathchardef\lne="7008
\mathchardef\sneq="7024
\mathchardef\spneq="7025
\mathchardef\sne="7028
\mathchardef\spne="7029
\mathchardef\ltms="706E
\mathchardef\tmsl="706F

\mathchardef\dbA="7041

%Euler Fraktur letters

%Capital roman double letters
\mathchardef\dbA="7041 
\mathchardef\dbB="7042 
\mathchardef\dbC="7043 
\mathchardef\dbD="7044 
\mathchardef\dbE="7045 
\mathchardef\dbF="7046 
\mathchardef\dbG="7047 
\mathchardef\dbH="7048 
\mathchardef\dbI="7049 
\mathchardef\dbJ="704A 
\mathchardef\dbK="704B 
\mathchardef\dbL="704C 
\mathchardef\dbM="704D 
\mathchardef\dbN="704E 
\mathchardef\dbO="704F 
\mathchardef\dbP="7050 
\mathchardef\dbQ="7051 
\mathchardef\dbR="7052 
\mathchardef\dbS="7053 
\mathchardef\dbT="7054 
\mathchardef\dbU="7055 
\mathchardef\dbV="7056 
\mathchardef\dbW="7057 
\mathchardef\dbX="7058 
\mathchardef\dbY="7059 
\mathchardef\dbZ="705A 

\def\nek{,\ldots,}
\def\sdp{\times \hskip -0.3em {\raise 0.3ex
\hbox{$\scriptscriptstyle |$}}} % semidirect product

%words in roman font

\def\dom{\mathop{\rm dom}\nolimits}

\def\min{\mathop{\rm min}}

%\def\rank{\rm rank}

%overlined math alphabet

%overlined Greek alphabet

%underlined math alphabet

\def\ut{{\underline t}}

%underline Greek alphabet
\def\ualp{{\underline\alp}}
\def\ubet{{\underline\bet}}

%math alphabet with hat

%Greek alphabet with hat

%roman with widetilde

\def\tilA{{\widetilde A}}

%Greek alphabet with widetilde

\def\ddownarrow{\big\downarrow \hskip-0.70em\raise
2pt\hbox {$\big\downarrow$}}
\def\longright #1#2 {\smash{\mathop{\hbox to
#1pt {\rightarrowfill}}\limits_{#2}}}
\def\sqr#1#2{{\vcenter{\hrule height.#2pt\hbox{\vrule
width.#2pt height#1pt \kern#1pt \vrule width.#2pt}
\hrule height.#2pt}}}

\def\buildrul#1\under#2{\mathrel{\mathop{\null#2}
\limits_{#1}}}

\def\boxit#1{\vbox{\hrule\hbox{\vrule\kern3pt
\vbox{\kern3pt#1 \kern3pt}\kern3pt\vrule}\hrule}}

\def\prodl{\prod\limits}

\def\subheading#1{\medskip\goodbreak\noindent{\bf
#1.}\quad}

\def\sect#1{\goodbreak\bigskip\centerline{\bf#1}
\medskip}
\def\pr{\smallskip\noindent{\bf Proof:\quad}}
\def\onumber #1{\ooalign{\hfil\raise.07ex\hbox{
\hfill$\scriptstyle \,#1$\hfil}
\cr\cr{$\bigcirc$}}}
\def\onumber c{\ooalign{\hfil\raise.07ex\hbox
{\hfill$\scriptstyle \,c$\hfil}
\cr\cr{$\bigcirc$}}}
\def\alpcirc {\ooalign{\hfil\raise.07ex
\hbox{\hfill$\scriptstyle\alp\;$\hfill}\cr\cr
{$\bigcirc$}}}

\def\longmapright #1#2 {\smash{\mathop{\hbox to
#1pt {\rightarrowfill}}\limits^{#2}}}
\def\longmapleft #1 #2 {\smash{\mathop{\hbox to
#1 pt {\leftarrowfill}}\limits^{#2}}}

\def\references#1{\goodbreak\bigskip\par\centerline
{\bf References}\medskip\parindent=#1pt}
\def\ref#1{\par\smallskip\hang\indent\llap{\hbox
to \parindent{#1\hfil\enspace}}\ignorespaces}

\def\back{{\raise 2.5pt\hbox{$\,\scriptscriptstyle
\backslash\,$}}}

\def\part{\partial}
\def\lwr #1{\lower 5pt\hbox{$#1$}\hskip -3pt}
\def\rse #1{\hskip -3pt\raise 5pt\hbox{$#1$}}
\def\lwrs #1{\lower 4pt\hbox{$\scriptstyle #1$}
\hskip -2pt}
\def\rses #1{\hskip -2pt\raise 3pt\hbox
{$\scriptstyle #1$}}

\def\<#1{\left\langle{#1}\right\rangle}

\def\subinbn{{\subset\hskip-8pt\raise 0.95pt
\hbox{$\scriptscriptstyle\subset$}}}

\def\llvdash{\mathop{\|\hskip-2pt
\raise 3pt\hbox{\vrule height 0.25pt width 1.5cm}}}

\def\lvdash{\mathop{|\hskip-2pt \raise 3pt\hbox
{\vrule height 0.25pt width 1.5cm}}}

\def\fakebold#1{\leavevmode\setbox0=\hbox{#1}%
  \kern-.025em\copy0 \kern-\wd0
  \kern .025em\copy0 \kern-\wd0
  \kern-.025em\raise.0333em\box0 }

\font\msxmten=msxm10
\font\msxmseven=msxm7
\font\msxmfive=msxm5
\newfam\myfam
\textfont\myfam=\msxmten
\scriptfont\myfam=\msxmseven
\scriptscriptfont\myfam=\msxmfive
\mathchardef\rhookupone="7016
\mathchardef\ldh="700D
\mathchardef\leg="7053
\mathchardef\ANG="705E
\mathchardef\lcu="7070
\mathchardef\rcu="7071
\mathchardef\leseq="7035
\mathchardef\qeeg="703D
\mathchardef\qeel="7036
\mathchardef\blackbox="7004
\mathchardef\bbx="7003
\mathchardef\simsucc="7025

\def\rhookup{{\fam=\myfam \rhookupone}}

\def\bigsquare{{\fam=\myfam\bbx}}

\font\tencaps=cmcsc10
\def\smallcaps{\tencaps}

\def\author#1{\bigskip\ce{\smallcaps #1}\medskip}

\def\upddots{\mathinner{\mkern
1mu\raise 1pt \hbox{.}\mkern 2mu \mkern
2mu \raise 4pt\hbox{.}\mkern 1mu \raise 7pt\vbox
{\kern 7 pt\hbox{.}}} }

\def\varchi{\ooalign{{\raise
1.385pt\hbox{$\chi$}}\crcr\hbox{--}\crcr}}

\def\trianarrow{{\raise 2pt\hbox to 0.50cm
{\hrulefill}\triangleright}}

%\input mssymb.tex

%${\vtop{\offinterlineskip\hbox
%{$r_1$}\hbox{$\sim$}}}$

\null
\overfullrule=0pt
\def\l{{\langle}}
\def\smallfrown{{\scriptscriptstyle \cap}}
\def\r{{\rangle}}
\def\C{{Cohen\ }}
\def\oi{{\overline i}}
\def\Suc{{\rm Suc}}

 \def\upp{{\buildrul \sim \under p}}
 \def\upp1{ {\buildrul \sim \under {p_1}}}
 \def\uppp{ {\buildrul \sim \under {p'_1}}}
 \def\up2{{ \buildrul \sim \under {p_2}}}
 
 \def\upn1{{ \buildrul \sim \under {p_{n1}}}}
\def\uq{{\buildrul \sim \under q}}
\def\uq1{{\buildrul \sim \under {q_1}}}
\def\ur1{{\buildrul \sim \under {r_1}}}
\def\uf{{\buildrul \sim \under f}}
\def\ug{{\buildrul \sim \under g}}
\def\ualp{{\buildrul \sim \under \alp}}
\def\ubet{{\buildrul \sim \under \bet}}

\def\ut{{\buildrul \sim \under t}}
\def\llvdash{\mathop{\|\hskip-1pt \raise 3pt\hbox{\vrule
height 0.25pt width 1cm}}}
\def\lvdash{\mathop{|\hskip-1pt \raise 3pt\hbox{\vrule
height 0.25pt width .5cm}}}
\def\ldash{\mathop{\|\hskip-1pt \raise 3pt\hbox{\vrule
height 0.25pt width .5cm}}}

\sect{ADDING A LOT OF COHEN REALS BY ADDING
A FEW}

\bigskip

\ce{Moti Gitik}
\medskip
\ce{School of Mathematical Sciences}
\ce{Sackler Faculty of Exact Sciences}
\ce{Tel Aviv University}
\ce{Tel Aviv 69978, Israel}

\vskip 1truecm
\dspace

A basic fact about Cohen reals is that
adding $\lam$ Cohen reals cannot produce more
than $\lam$ of Cohen reals.  More
precisely, if $\l r_\alp | \alp < \lam \r$
are $\lam$-Cohen generic reals over $V$,
then in $V[\l r_\alp | \alp < \lam \r ]$
there is no $\lam^+$-Cohen generic real
over $V$.

But if instead of dealing with one universe
$V$ we consider two, then the above may no
longer be true.

The purpose of this paper is to produce
models $V_1 \subseteq V_2$ such that adding
$\kap$-many Cohen reals to $V_2$ adds
$\lam$-Cohen to $V_1$, for some $\kap <
\lam$.  We deal mainly with the case when
$V_1, V_2$ have same cardinals and satisfy
GCH.

Let us state the principal results:
\proclaim Theorem 1.1.  Suppose that $V$
satisfies GCH, $\kap = \bigcup\limits_{n <
\ome} \kap_n$ and $\bigcup\limits_{n <
\ome} o(\kap_n) = \kap$ (where $o(\kap_n)$
is the Mitchell order of $\kap_n$).  Then
there exists a cardinal preserving generic
extension $V_1$ of $V$ satisfying GCH and
having the same reals as $V$ does, so that
adding $\kap$-many Cohen reals over $V_1$
produces $\kap^+$-many Cohen reals over
$V$.

\proclaim Theorem 2.1.  Suppose that $V
\models GCH$.  Then there exists a
cofinality preserving extension $V_1 \models
GCH$ so that adding a Cohen real to $V_1$
produces $\aleph_1$-Cohen reals over $V$.

\proclaim Theorem 2.2.  There is a pair
$(W, W_1)$ of a generic cofinality
preserving extensions of $L$ such that $W
\subseteq W_1$ and $W_1$ contains a perfect
set of $W$-reals which is not in $W$.

It is a slight improvement of a result from
Velickovic-Woodin [Ve-Wo] since there
$\aleph_2$ was collapsed.

Some large cardinal assumptions are needed
for 1.1.  Since, if every countable subset
of $\kap$ of $V_1$ can be covered by
a countable set of $V$ and $V, V_1$ have
same reals then it is impossible to produce
$\kap^+$-Cohen reals over $V$ by adding
$\kap$-Cohen reals over $V_1$.

In Section 1, we use  
the B. Velickovic and H. Woodin [Ve-Wo] idea of adding a club
avoiding ordinals which have countable
cofinality in $V$.  

The present work was inspired by a work by 
Velickovic-Woodin [Ve-Wo] and a question by
Fremlin [Fr].  We are grateful to B.
Velickovic for discussions on the subject.

\sect{1.~~Models with the same reals}

\proclaim Theorem 1.1.  Suppose that $V$
satisfies GCH, $\kap = \bigcup\limits_{n <
\ome} \kap_n$ and $\bigcup\limits_{n <
\ome} o(\kap_n) = \kap$ ($o(\kap_n)$ is the
Mitchell order of $\kap_n$).  Then there
exists a cardinal preserving generic
extension $V_1$ of $V$ satisfying GCH and
having the same reals as $V$ does so that
adding $\kap$-many Cohen reals over $V_1$
produces $\kap^+$ many Cohen reals over
$V$.

\pr Rearranging the sequence $\l \kap_n | n
< \ome\r$, we can assume that
$o(\kap_{n+1})\ge \kap_n^+$.  First let us
change the cofinality of $\kap_{n + 1}$ to
$\kap_n^+$ for every $n, 0 < n < \ome$
using the Magidor forcing [Ma].  It does
not add new subsets to $\kap_0$ and
preserves GCH.  The second step will be to
add a club $C \subseteq \kap^+$ such that
for every $\alp \in C$ ($cf \alp)^V >
\ome$.  We use the usual forcing  
$\calP = \{ c | c \subseteq \kap^+, \ |c|
< \kap^+, \ c$ is closed and $c$ consists
of ordinals $\alp$ such that $(cf
\alp)^V$ is singular$\}$ ordered by
end extensions.

\proclaim Claim 1.1.1.  The forcing $\calP$
over $V^{\hbox{Changing of Cofinalities}}$
is $(\kap, \infty)$-distributive.

\pr Let $\del < \kap$.  In order to apply
the standard argument, we need to find for every
regular $\del < \kap$  an elementary
submodel $N$, $|N| = \del$ such that $cf^
V(\sup (N \cap \kap^+))$ is a singular
cardinal of cofinality $\del$.  For this it
is enough to show that the set of $V$
$S^\xi_{\kap^+} = \{ \alp < \kap^+ | cf^V \alp
= \xi\}$
remains stationary after changing
cofinalities forcing for every regular in
$V$ cardinal $\xi < \kap$.  But notice
that, since $cf^V \kap = \ome$, the
changing cofinalities forcing satisfies
$\kap^+$-c.c. (it is even
$\kap$-centered).  So each $S^\xi_{\kap^+}$ remains
stationary and we are done.\hb$\bigsquare$

Notice, that by Claim 1.1.1, $\kap^+$
remains a cardinal, since otherwise it
would change the cofinality to some $\del <
\kap$ and so $\calP$ would add a new function
from $\del$ to $\kap^+$ which is
impossible.

Let $C$ be a generic club of $\kap^+$ added
by $C$.  Set $V_1 = V$ [changing
cofinalities, $C$].  The following is
obvious.

\proclaim Claim 1.1.2.  For every countable
in $V$ set $t$ and $\alp \in C$, $t\cap
\alp$ is
bounded in $\alp$.

Fix $\l f_\alp | \alp < \kap^+\r \in V$ an
increasing (mod finite) sequence of
functions in $\prodl_{n < \ome} (\kap_{n +
1}^+ \backslash \kap_{n + 1})$.

Now we force over $V_1$ $\kap$ \C reals $\l
r_i | i < \kap \r$.  Our aim is to
construct from them another sequence $\l
s_\alp | \alp < \kap^+ \r$ which will be a
sequence of $\kap^+$-\C reals over $V$.

Let us first split $\l r_i | i < \kap \r$
into two sequences of length $\kap$
denoted by $\l r_i | i < \kap\r$ and
$\l r'_i | i < \kap\r$. 

Let $\alp < \kap^+$.  We define $s_\alp$ as
follows
\subheading{Case 1}  $\alp \in C$.\hb
For every $n < \ome$ set
$$s_\alp (n) = r_{f_\alp (n)} (0)\ .$$

\subheading{Case 2} $\alp \not\in C$.\hb
Let $\alp^*$, $\alp^{**}$ be two successor
points of $C$ so that $\alp^* < \alp <
\alp^{**}$.  Let $\l \alp_i | i < \kap \r$
be some fixed from the beginning
enumeration of the interval $(\alp^*,
\alp^{**})$.  Then for some $\oi < \kap$
$\alp = \alp_\oi$.  We consider the \C real
$r'_\oi$.  Let $k (\oi)$ be the least $k <
\ome$ so that $r'_\oi (k) = 1$. Set
$$s_\alp (n) = r_{f_\alp (k (\oi) + n)}
(0)$$
for every $n < \ome$.

\proclaim Claim 1.1.3.  $\l s_\alp | \alp <
\kap^+ \r$ is a sequence of $\kap^+$ \C
reals over $V$.

\pr Let $C(\lam \times \ome, 2)$ denote the \C
forcing for adding $\lam$ \C reals, i.e.
$C(\lam \times \ome, 2) = \{ p \mid p \
\hbox{is a finite function from}\ \lam
\times \ome\ {\rm into}\ 2\}$.

Suppose that $I \in V$ is a maximal
antichain in $C(\kap^+ \times \ome, 2)$.
Let $ t(I) = \{ \alp < \kap^+| \alp$ appears
in some element of $I\}$.  By c.c.c.
$t(I)$ is a countable subset of $\kap^+$.

We need to show that every condition
$$\l p, p'\r \in C(\kap \times \ome, 2)
\times C(\kap \times \ome, 2)$$
can be extended to some $\l q, q'\r$
forcing
$$''\l 
\vtop{\offinterlineskip\hbox{$s_\alp$ } \hbox{$\scriptstyle\sim$}}
\mid
\alp \in \check t (I) \r\ \hbox{extends
some element of}\ \check I''\ .$$

Let us assume for simplicity that $\l p,
p'\r = \emptyset$.  Recall also that we
split $\kap$ \C reals into two $\l r_\alp |
\alp < \kap \r$ and $\l r'_\alp | \alp <
\kap\r$.

By Claim 1.1.2, $t(I)$ can have a nonempty
intersection only with finitely many
intervals $[\alp^*, \alp^{**})$ where
$\alp^*$, $\alp^{**} \in C$ and $\alp^{**}
= \min (C \backslash(\alp^* + 1))$. Suppose
for simplicity that there are only two such
intervals $[\alp^*_1, \alp_1^{**})$ and
$[\alp_2^*, \alp_2^{**})$, where $\alp_1^*
< \alp_2^*$.  First of all let us deal with
$s_{\alp_1^*}$, $s_{\alp_2^*}$ (or with
names of these two reals).  Recall that
$s_\alp (n) = r_{f_\alp(n)} (0)$ for every
$\alp \in C$ and $n < \ome$.  Pick some
$n^* < \ome$ such that for every $n \ge
n^*$ $f_{\alp_1^*} (n) < f_{\alp_2^*} (n)$.
Now let $p \in C(\kap \times \ome, 2)$ be a
condition such that dom $p = \{ (\bet, 0) \in
\kap \times \ome|$ for some $n < n^*$
$f_{\alp_1^*} (n) = \bet$ or $f_{\alp_2^*}
(n) = \bet\}$.  We define $s_{\alp_1^*}
\rhookup n^*$, $s_{\alp_2^*} \rhookup n^*$
according to $p$.  Intuitively, we decide
the values of $s_{\alp_1^*}$,
$s_{\alp_2^*}$ over $n$'s where they may
disturb one another.  Now find some $b \in
I$ stronger than $\l s_{\alp_1^*} \rhookup
n^*, s_{\alp_2^*} \rhookup n^*\r$.  Let us
extend $p$ to a condition $q$ forcing 
$b$ to be satisfactory, i.e.
$$q \ldash \l 
\vtop{\offinterlineskip\hbox{$s_{\alp_1^*}$ } \hbox{$\scriptstyle\sim$}}
, 
\vtop{\offinterlineskip\hbox{$s_{\alp_2^*}$ } \hbox{$\scriptstyle\sim$}}
\r\ {\rm
 extends}\ \l \check
b(\alp_1^*), \check b(\alp_2^*) \r\ .$$

This is not problematic since for $n \ge
n^*$ $f_{\alp_1^*} (n) \ne f_{\alp_2^*}
(n)$ and also $f_{\alp_j^*} (n) \in \kap_{n
+ 1}^+ \backslash \kap_{n + 1}$ $(j \in \{
1,2\})$ hence $f_{\alp_1^*} (n) \not\in \{
f_{\alp_1^*} (m), f_{\alp_2^*} (m) | m <
n\}$.

Now we need to take care of the rest
coordinates of $b$.  Let $\{ \alp_{j1}\nek
\alp_{jk_j}\}$ be the increasing
enumeration of the coordinates of $b$ in
the interval $(\alp^*_j, \alp_j^{**})$ for
$j \in \{ 1,2\}$ and $k_j < \ome$.  For $j
\in \{1,2\}$ and $\ell < k_j$ let
$\oi_{j,\ell} < \kap$ be the index of
$\alp_{j\ell}$ in the enumeration of the
interval $(\alp_j^*, \alp_j^{**})$ used in
Case 2 of the definition.  Recall that
$s_{\alp_{j\ell}} (n) =
r_{f_{\alp_{j\ell}} (k (j,\ell) + n)}
(0)$, where $n < \ome$ and $k(j, \ell)$ is
the least $k < \ome$ such that
$r'_{\oi_{j,\ell}} (k) = 1$.

In order to prevent collisions let us pick
$m^* < \ome$ big enough so that for every
$n \ge m^*$, $j \in \{ 1,2\}$, $\ell_j,
\ell'_j < k_j$, $\ell_j < \ell'_j$
$$f_{\alp_1^*} (n) < f_{\alp_{1 \ell'_1}}
(n) < f_{\alp_2^*} (n) < f_{\alp_{2\ell_2}}
(n)$$
and 
$$f_{\alp_{j \ell_j}} (n) < f_{\alp_{j
\ell'_j}} (n)\ .$$

There is no problem in  finding such $m^*$
since there are conditions on finitely
many functions only. 

Now for every $j \in \{ 1,2\}$ and $\ell <
k_j$ define the first $m^*$ values of
$r'_{\oi_{j,\ell}}$ to be 0.  Hence
$k(j,\ell)$ will be forced  to be $\ge
m^*$.  But then the reals $r_{f_\alp
(k(j,\ell) + n)}$ will be all different,
where $\alp \in t(I)$.  Translate this into
a
condition in $C(\kap \times \ome, 2) \times
C(\kap\times \ome, 2)$ and we are done.\hb
$\bigsquare$ of the claim.\hb
$\bigsquare$

The same result is true for random reals.
The proof is similar but a little more
involved, since random reals in contrast to
Cohen reals may depend on $\ome$-many coordinates
instead of finitely many in the Cohen case.

\proclaim Theorem 1.2.  Suppose that $V$
satisfies GCH, $\kap = \bigcup\limits_{n <
\ome} \kap_n$ and $\bigcup\limits_{n <
\ome} o(\kap_n) = \kap$.  Then there exists
a cardinal preserving generic extension
$V_1$ of $V$ satisfying GCH and having the
same reals as $V$ does so that adding
$\kap$-many random reals over $V_1$
produces $\kap^+$-many random reals over
$V$.
 
\subheading{Sketch of the proof}

We define $\l s_\alp | \alp < \kap^+\r$ as
in 1.1.  The proof proceeds as in 1.1 up to
the point of taking care of coordinates of
$b$ which differ from $\alp_1^*, \alp_2^*,
\alp_1^{**}, \alp_2^{**}$.  Here there may
be
infinitely many coordinates of $b$ inside
intervals $(\alp^*_1, \alp_1^{**})$ and
$(\alp_2^*, \alp_2^{**})$.  Let us analyze
the situation more closely and figure out
reasons for possible collisions.  Let
$\{\alp_{j\ell} | \ell < k_j \le \ome\}$ be
an enumeration of coordinates of $b$ 
in the interval $(\alp_j^*, \alp_j^{**})$
for $j = 1,2$.  For $j \in \{ 1,2\}$ and
$\ell < k_j$ let $\oi_{j,\ell} < \kap$ be
the index of $\alp_{j\ell}$ in the
enumeration of the interval $(\alp_j^*,
\alp_j^{**})$ used in Case 2 of the
definition of $\l s_\alp | \alp <
\kap^+\r$. Recall that
$$s_{\alp_{j\ell}} (n) = r_{f_{\alp_{j\ell}}
(k(j,\ell) + n)} (0)\ ,$$
where $n < \ome$ and $k(j,\ell)$ is the
least $k < \ome$ such that
$r'_{\oi_{j,\ell}} (k) = 1$.  For every
$j_1, j_2 \in \{ 1,2\}$, $\ell_1 < k_{j_1}$,
$\ell_2 < k_{j_2}$ and $m,n < \ome$, let
$$c(j_1, j_2, \ell_1, \ell_2, m) = \|
s_{\alp_{j_1\ell_1}} (n) \ne s_{\alp_{j_2
\ell_2}} (m)\|\ .$$
Then there are only finitely many $(j_1, j_2,
\ell_1, \ell_2, n_1, n_2)$ such that $b \le
c(j_1, j_2, \ell_1, \ell_2, n_1, n_2)$.
Now as in 1.1 we pick $m^* < \ome$ large
enough  to insure that above it all
relevant functions are different.  Finally,
for every $n \le m^*$ we set $r'_{\oi_{j,
\ell}} (n) = 0$ where $(j, \ell)$ appears in
the finite set above.\hb
$\bigsquare$

Let us show that some large cardinals are
needed for the previous results.  First the
following trivial observations:

\proclaim Proposition 1.3.  Suppose that $V_1
\supseteq V$, $V_1$ and $V$ have the same
reals and $([\lam]^{\le \ome} )^V$ is
unbounded in $([\lam]^{\le \ome} )^{V_1}$
for some $\lam$.  Then $([\lam]^{\le
\ome})^V = ([\lam]^{\le \ome} )^{V_1}$.

\proclaim Proposition 1.3.1.  Assume that
$V_1 \supseteq V$, $V_1$ and $V$ have the
same reals. Suppose that for some uncountable
cardinal $\kap$ of $V_1$ adding $\kap$-\C
(random) reals to $V_1$ adds
$(\kap^+)^{V_1}$-\C (random) reals to $V$,
then $([\kap^+]^{\le \ome_1} )^V$ is
bounded in $([\kap^+]^{\le \ome})^{V_1}$.

\pr Notice that $\ome_1^V = \ome_1^{V_1}$.
If $A$ is a countable subset of $\kap^+$ in
$V_1$ and $A^* \in V$, $A^* \supseteq A$
and $|A^*|^V = \aleph_1$, then for some
$A^{**} \subseteq A^*$, $A^{**} \in V$ and
$A^{**}$ countable $A^{**} \supseteq A$.
Now the claim follows by 1.3 and c.c.c. of
the forcing.\hb
$\bigsquare$

By the Covering Lemma of Dodd-Jensen, we
now obtain the following

\proclaim Corollary 1.3.2.  Under the
assumptions of 1.3.1 and if $\aleph_2^V =
\aleph_2^{V_1}$ then there is an inner
model with a measurable cardinal.

We do not know the exact strength. 

\subheading{Question}  Is
$\bigcup\limits_{n < \ome} o(\kap_n) =
\kap$ really needed for the results 
1.1, 1.2?
\hb
 If we relax our assumptions and allow 
 at least cardinals $\aleph_1$ and
$\aleph_2$  to collapse then no large cardinal
assumptions are needed.  Namely the
following holds.

\proclaim Theorem 1.4.1.  Suppose $V$ is a
model of GCH then there is a generic
extension $V_1$ of $V$ such that all the
cardinals except $\aleph_1$ are preserved,
GCH holds and adding of $\aleph_\ome$-Cohen
(random) reals over $V_1$ produces
$\aleph_{\ome+1}$ many of them over $V$.

\pr We need a club $C \subseteq
\aleph_{\ome + 1}$ satisfying Claim 1.1.2
of Theorem 1.1.  If we have it then the
rest of the poof is as in this theorem.
The construction is as in Velicovich-Woodin
[Ve-Wo].  First collapse  $\aleph_1$ to
$\aleph_0$.  Then add a club $C_1$ to new
$\aleph_1$ (former $\aleph_2$) consisting
of points of cofinality $\aleph_1$ in $V$.
Over $\aleph_2$  force a club $C_2$
such that every $\alp \in C_2$ of countable
cofinality has cofinality $\aleph_1$ in
$V$.  Continue in the same fashion over all
$\aleph_n$'s.  Finally add $C \subseteq
\aleph_{\ome+1}$ as desired.\hb
$\bigsquare$

If one wants to keep the same reals, then
just replace the collapse of $\aleph_1$ by the
Namba forcing.

By the same lines but using stronger
initial assumptions adding $\kap$ many reals
\C or random may produce $\lam$ many of
them for $\lam$ much larger than $\kap^+$.

\proclaim Theorem 1.5.  Suppose $\kap$ is a
strong cardinal, $\lam \ge \kap$ regular and
GCH.  Then there exists a cardinal
preserving extension $V_1$ having no new
reals such that adding $\kap$ many \C
(random) reals over $V_1$ produces $\lam$
many of then over $V$.

\pr Using Radin forcing, add a club $C_\kap
\subseteq \kap$  consisting of regular
cardinals such that $\kap$ is still strong
in $V[C_\kap]$.  Now, as in 1.1, force a
club $C_{\kap^+}$ consisting of $\alp$'s
below $\kap^+$ such that $(cf \alp)^V$ is
singular or equal to $\kap$.

Then we force a club $C_{\kap^{++}}
\subseteq \kap^{++}$.  It will consist of
$\alp$'s so that $(cf \alp)^V$ is singular
or $cf \alp = \kap$ or $cf \alp = \kap^+$.
As in 1.1 this forcing does not add new
subsets to $\kap^+$.  Now we continue this
iteration all the way to $\lam$.  The full
support is used at limit stages.  There is
no problem with distributivity since we
can always  take elementary submodels of
regular cardinalities $\ge \kap$.

Finally we will have a club $C_\lam$ of
$\lam$ consisting of ordinals $\alp < \lam$
such that $(cf \alp)^V$ is singular or $cf
\alp$ is a regular cardinal $\ge \kap$.  No
new subsets of $\kap$ are added by such
extension.  Let us denote the resulting
model by $V'$.  There was an extender of
length $\lam$ over $\kap$ in $V$.  Use it
as in [Gi-Ma1,2] to blow the power of
$\kap$ to $\lam$ and simultaneously to change
its cofinality to $\aleph_0$. Let $V_1$ be
such an extension.  By [Gi-Ma1,2] there is a
sequence $\l \kap_n | n < \ome\r$ of
regular cardinals cofinal in $\kap$ and a
sequence of functions $\l  f_\alp | \alp <
\lam \r$ in $\prodl_{n < \ome} (\kap_{n +
1} \backslash \kap_n^+)$ such that
$f_\alp(n) < f_{\alp'} (n)$ for all but
finitely many $n$'s whenever $\alp <
\alp'$.  Also if $\kap < \tau \le \lam$,
$\tau$ is 
regular, $\alp \in C_\tau$ and $cf \alp =
\aleph_0$ (in $V_1$), then $\alp$ is regular
in $V$.  Hence Claim 1.2 holds for
$C_\tau$. 

We now construct  $\lam$-Cohen (random)
reals as in Theorem 1.1 (1.2) using
$C_\lam$ and $\l f_\alp | \alp < \lam \r$.

Cases 2 of the definition of $\l s_\alp |
\alp < \lam\r$ is now problematic since the
cardinality of an interval $(\alp^*,
\alp^{**})$ (using the notation of 1.1) may
now be above $\kap$ but we have only $\kap$
many Cohen reals to play with.  Let us
proceed as follows in order to overcome
this.

We force over $V_1$ $\kap$-\C reals as
before.  But
let us rearrange them as $\l r_{n, \alp} |
\alp < \kap, n < \ome\r$, $\l r_\eta | \eta
\in [\kap]^{< \ome} \r$.  We define by
induction on levels a tree $T \subseteq
[\lam]^{< \ome}$, its projection $\pi (T)
\subseteq [\kap]^{< \ome}$ and for $n <
\ome$ and $\xi \in {\rm Lev}_n T$ a real
$s_\xi$.  The union of the levels of $T$
will be $\lam$, so $\l s_\xi | \xi <
\lam\r$ will be defined.

Set ${\rm Lev}_0 T = \l\ \r = {\rm Lev}_0
\pi(T)$.  Define Lev$_1 T = C_\lam$ and
Lev$_1\pi (T) = \{0\}$, i.e. $\pi (\l \alp
\r ) = \l 0 \r$ for every $\alp \in
C_\lam$.  For $\alp \in C_\lam$ we define
$$s_\alp (m) = r_{1, f_\xi (m)} (0)$$
for every $m < \ome$.

Suppose now that $n > 1$ and $T \rhookup
n$, $\pi (T \rhookup n)$ are defined.
Let $\eta \in T \rhookup n - 1$, $\alp^*,
\alp^{**} \in \Suc_T (\eta)$ and $\alp^{**}
= \min (\Suc_T (\eta) \backslash \alp^*)$.
Define $\Suc_T (\eta ^\smallfrown \alp^{**})$ if it
is not yet defined.

\subheading{Case A} $|\alp^{**} \backslash
\alp^*| \le \kap$.\hb
Then set $\Suc_T (\eta ^\smallfrown \alp^{**}) =
\alp^{**} \backslash \alp^*$ and $\Suc_T
(\eta ^\smallfrown \alp^{**\smallfrown}
\alp) = \emptyset$.  Define $\Suc_{\pi(T)}
(\pi (\eta ^\smallfrown \alp^{**})) =
|\alp^{**} \backslash \alp^*|$ and $\pi
(\eta ^\smallfrown \alp^{**}) ^\smallfrown i
= \emptyset$ for every $\alp \in \alp^{**}
\backslash \alp^*$, and $i < | \alp^{**}
\backslash \alp^* |$.  Let us define $s_\alp$
for $\alp \in \alp^{**} \backslash \alp^*$.
Fix some enumeration $\l \alp_i | i < \rho
\le \kap\r$ of $\alp^{**} \backslash
\alp^*$.  Let $\oi$ be so that $\alp_\oi =
\alp$. Set $s_\alp(m) = r_{n,f_\alp (k +
m)} (0)$ for every $m < \ome$, where $k$ is
the least $\ell$ such that $r'_{\pi (n
^\smallfrown \alp^{**}) ^\smallfrown \oi} (k)
\ne 0$.

\subheading{Case B} $|\alp^{**} \backslash
\alp^*| > \kap$ and $cf \alp^{**} < \kap$.
\hb
Let $\rho = cf \alp^{**}$.  Pick an
increasing continuous cofinal in $\alp^{**}$
sequence $\l \alp^{**}_\nu | \nu < \rho\r$
with $\alp_0^{**} > \alp^*$.  Set $\Suc_T
(\eta ^\smallfrown \alp^{**}) = \{
\alp_\nu^{**} | \nu < \rho\}$ and $\Suc_{\pi(T)}
(\pi
(\eta ^\smallfrown \alp^{**})) = \{ \nu | \nu
< \rho \} = \rho$.  Set $s_{\alp^{**}_0}
(m) = r_{n, f_{\alp^{**}_\nu}} (k + m) (0)$
for every $m < \ome$ and $\nu < \rho$,
where $k$ is the least such that $r'_{\pi
(\eta ^\smallfrown \alp^{**}) ^\smallfrown
\nu} (k) \ne 0$.

\subheading{Case C} $cf \alp^{**} >
\kap$.\hb
Let $\rho$ and $\l \alp_\nu^{**} | \nu <
\rho\r$ be as in Case B.  Set $\Suc_T (\eta
^\smallfrown \alp^{**}) = \{ \alp_\nu^{**} |
\nu \in C_\rho\}$ and $\pi (\eta
^\smallfrown \alp^{**}) = \{ 0\}$.  For
every $\nu \in C_\rho$ and $m < \ome$
define $s_{\alp_\nu^{**}} (m) =
r_{n,f_{\alp_\nu^{**}}} (k + m) (0)$ where
$k$ is the least such that $r'_{\pi (\eta
^\smallfrown \alp^{**})^\smallfrown 0} (k)
\ne 0$.

By the definition $T$ is a well-founded tree.
Clearly, every $\alp < \lam$ appears in $T$
at some level.

\proclaim Claim 1.5.1.  $\l s_\alp | \alp <
\lam\r$ are $\lam$-\C reals over $V$.

\pr Let $I, t(I)$ be as in 1.1.1.  If $n_1
\ne n_2 < \ome$, then there are no
collisions between $t(I) \cap {\rm
Lev}_{n_1} (T)$ and $t(I) \cap {\rm
Lev}_{n_2} (T)$ since different \C reals
are used over different levels.

Suppose that for some $n < \ome$, $t(I) \cap
{\rm Lev}_n (T)$ is infinite.  Then find
the least level $m \le n$ such that for
infinitely $\alp \in {\rm Lev}_m (T)$
$\alp$ has a successor on the level $n$ which
belongs to $t(I)$.  Or in other setting,
there are infinitely many $\alp_k^* <
\alp_k^{**}$ immediate successors in the
ordering of   the
$m$-th level of $T$ such that $(\alp^*_k,
\alp^{**} _k) \cap t(I) \ne \emptyset$.  By
Claim 1.1.1, $m > 1$.  But then there is
$\alp \in {\rm Lev}_{m-1} (T)$ such that for
infinitely many $k$'s $\alp^{**}_k \in
\Suc_T (\alp)$.  
So we are in a situation of Case A or Case
B, since by Claim 1.1.1, $cf \alp$ cannot be
$> \kap$.  Let $k_1 \ne k_2$ be two such
$k$'s.  Then, by the definition of $\pi$,
$\pi (\eta^\smallfrown \alp ^\smallfrown
\alp^{**}_{k_1}) \ne \pi (\eta ^\smallfrown
\alp ^\smallfrown \alp^{**}_{k_2})$, where
$\eta \in T \rhookup m-1$ is the set
of predecessors of $\alp$.  This implies that
also over level $n$ different $r'_\nu$'s
are used to define $s_\alp$'s $(\alp \in t(I)$)
for all but finitely many $\alp$'s in $t(I)$.
The rest of the argument is as in 1.1.2.\hb
$\bigsquare$

If we allow many cardinals
between $V$ and $V_1$ to collapse, then using [Gi-Ma1,
Sec. 2] one can obtain the following:

\proclaim Theorem 1.6.  Suppose that there
is a strong cardinal and GCH holds.  Let
$\alp < \ome_1$, then there is a model $V_1
\supseteq V$ having the same reals as $V$
satisfying GCH below $(\aleph_\ome)^{V_1}$
such that adding $(\aleph_\ome)^{V_1}$-many
\C reals to $V_1$ produces $(\aleph_{\alp +
1})^{V_1}$-many of them over $V$.

The proof is straightforward using 1.5 and
[Gi-Ma1].  Notice that in this model $V_1$
$\aleph_\ome$ is much bigger than
$(\aleph_\ome)^V$.  So the following is
natural.

\subheading{Question}  Is it possible to
have $V_1 \supseteq V$ so that
\item{(a)} $V_1$ and $V$ have the same
reals,
\item{(b)} $V_1$ and $V$ have the same
cardinals,

\item{(c)} adding $\aleph_\ome$-\C reals to
$V_1$ adds more than $\aleph_\ome$-\C
reals to $V$?

By 1.3 it requires large cardinals.  A
weaker version which may be consistent with
just
 ZFC:

\subheading{Question}  Is it possible to
have $V_1 \supseteq V$ so that (b) and (c)
of the previous question hold?

\sect{2.~~Models with the same cofinality
function but different reals}

In this section we are going to show the
following:

\proclaim Theorem 2.1.  Suppose that $V
\models $GCH.  Then there exists a
cofinality preserving extension $V_1
\models $GCH so that adding a Cohen (or
random) real to $V_1$ produces $\aleph_1$
\C (or random) reals over $V$.

A similar construction will be used to give
another proof of a result of B. Velickovic
and H. Woodin related to an old open
question of K. Prikry: whether $R^L$ could
contain a perfect set and not to be equal
to the set of all reals.

\proclaim Theorem {\rm [Ve-Wo]}.  There is
a pair $(W, W_1)$ of generic extensions of
$L$ such that $W \subseteq W_1$,
$\aleph_1^W = \aleph_1^{W_1}$ and $W_1$
contains a perfect set of $W$-reals but not
all the reals of $W_1$ are in $W$.

$\aleph_2$ is collapsed in their
construction.

We will get  a slightly stronger result:

\proclaim Theorem 2.2.  There is a pair
$(W, W_1)$ of generic cofinality preserving
extensions of $L$ such that $W \subseteq
W_1$, $W_1$ contains a perfect set of
$W$-reals and not all the reals of $W_1$
are in $W$.

The basic idea of the proofs will be to
split $\ome_1$ in 2.1, 2.2 into $\ome$ sets
such that none of them will contain an
infinite set of $V$.  Then something like
in Section 1 will be used for producing \C
reals.

It turned out however that just not
containing an infinity set of $V$ is not
enough.  We will use a stronger property.
As a result the forcing turns out to be more
complicated.

We are going to define the forcing
sufficient for proving Theorems 2.1, 2.2.

Fix a nonprincipal ultrafilter $U$ over
$\ome$.

\subheading{Definition 2.1.1}  Let $\calP_U$
be the Prikry (or in this contest Mathias)
forcing with $U$, i.e.
$$\l p, A\r \in \calP_U \quad {\rm iff}
\quad p \in [\ome]^{<\ome}\ ,$$
$A \in U$ and $\sup p < \min A$. $\l p, A \r
\le \l q, B\r$ iff $q$ is an end extension
of $p,q \backslash p \subseteq A$ and $B
\subseteq A$.

Define also $\l p, A\r \le^* \l q, B\r$ iff
$p = q$ and $A \supseteq B$.  We call
$\le^*$ a direct or $*$-extension.

The following are the basic facts on this
forcing that will be used further.

\subheading{Fact 1}  The generic object for
$\calP_U$ is generated by a real.

\subheading{Fact 2}
$\l \calP, \le \r$ satisfy c.c.c.

\subheading{Fact 3} If $\l p, A\r \in
\calP_U$ and $b$ is a finite subset of
$\ome \backslash \max p$, then a
$*$-extension $\l p, A \backslash b\r$ of
$\l p, A\r$ forces the  generic real to be
disjoint to $b$.

Now let us turn to the definition of forcing
that is going to split (or color) $\ome_1$
into $\ome$ pieces (or colors) without
infinite homogeneous set in $V$.

\subheading{Definition 2.1.2}

$p = \l p_0,
%\vtop{\offinterlineskip\hbox{$p_1$ } \hbox{$\scriptstyle\sim$}}
\upp1 \r \in \calP$ iff
 \item{(1)} $p_0 \in \calP_U$
\item{(2)} $\upp1$ is a $\calP_U$ name of a
coloring of some ordinal $\alp < \ome_1$
into $\ome$ colors such that the following
holds:

\itemitem{\rm (2a)} for every $\bet < \alp$
$\upp1(\bet)$ is a partial function from
$[\ome]^{<\ome}$ into $\ome$ and if $\bet_1
\ne \bet_2 < \alp$ then $rng({\buildrul
\sim \under {p_1}} (\bet_1))
\cap rng(\upp1 (\bet_2))$ is finite.

\itemitem{\rm (2b)}(*) for every countable $I
\subseteq \alp$, $I \in V$, $p'_0 \ge p_0$
and a finite $J \subseteq \ome$ there
exists a finite $a(I,J, p'_0) \subseteq
\alp$ such that for every finite $b
\subseteq I \backslash a(I, J, p'_0)$ there
is $p_0^{''} \ge^* p'_0$
$p_0^{''} {\buildrul \calP_U \under
\llvdash} \big( \forall k \in J\ \forall
\bet \in b (\upp1 (\bet) \ne k) \big)$ and
$(\forall \bet_1 \ne \bet_2 \in b\ \upp1
(\bet_1) \ne \upp1 (\bet_2))$, \hb
i.e.
$p_0^{''}$ prevents any  element of $b$ from
being colored by the colors in $J$ and no two
elements of $b$ are colored by the same
color.

Let us further denote $\alp$ by $\ell th$
$(p)$ or $\ell th$ $(\upp1)$ (the length of
$\upp1$).  Notice that such $\alp$ is not a
$\calP_U$-name.  $p_0$ already decides its
value.  Let us also denote by $\l
\vtop{\offinterlineskip\hbox{$I_{p,n}$ } \hbox{$\scriptstyle\sim$}}
| n < \ome
\r$ the splitting of $\alp$ induced by the
coloring $\upp1$.

The condition (2b)$(*)$ appears technical but it
will be crucial for producing numerous \C
reals.

\subheading{Definition 2.1.3}  Let $p = \l
p_0, \upp1\r$, $q = \l q_0, \uq1 \r \in
\calP$.  Define $\l p_0, \upp1\r \le \l
q_0, \uq1\r$ iff $p_0 \le q_0$, $\ell th$
$(p) \le \ell th (q)$ and $q_0 {\buildrul
\calP_U \under \llvdash}^{''}$ for every $n
< \ome$, 
$
\vtop{\offinterlineskip\hbox{$I_{p,n}$ } \hbox{$\scriptstyle\sim$}}
 = 
\vtop{\offinterlineskip\hbox{$I_{q,n}$ } \hbox{$\scriptstyle\sim$}}
 \cap \ell th (p)^{''}$.
Define also $\l p_0, \upp1\r \le^* \l q_0,
\uq1\r$ iff $p_0 \le^* q_0$ and $\l p_0,
\upp1 \r \le \l q_0, \uq1\r$.

\proclaim Lemma 2.1.4.  Suppose that $\l
p_0, \upp1\r \ldash^{''} \ualp$ is an
ordinal". Then there are a $\calP_U$-names
$\ubet$ and $\uppp$ such that $\l p_0,
\upp1\r \le^* \l p_0, \uppp\r$ and $\l p_0,
\uppp\r \ldash '' \ualp = \ubet''$.

\subheading{Remark}  In the situation of
the usual iteration this lemma is trivial.
But here we have in addition to satisfy the
condition (2b)$(*)$ which complicates 
matters.

\pr Suppose for simplicity that $\l p_0,
\upp1\r$ is just the weakest condition,
i.e. $\l\l\l\ \r, \ome \r, \emptyset \r$.  Let
$\l N_n | n < \ome\r$ be an increasing
sequence of countable elementary submodels
of a large enough portion of the universe
such that
\item{(a)} $\calP, \ualp \in N_0$,
\item{(b)} $N_n \in N_{n + 1}$ for every $n
< \ome$.
\hb
Set $N = \bigcup\limits_{n < \ome} N_n$.
Let $\del_n = N_n \cap \ome_1$ and $\del =
\bigcup\limits_{n < \ome} \del_n = N \cap
\ome_1$.

Fix an enumeration
$\l s_n | n < \ome \r \in N_0$ of
all increasing finite sequences of natural
numbers such that $n \le m$ implies max
$s_n \le \max s_m$.
Fix also a partition $\l J_{\vec n} | \vec
n \in [\ome]^{< \ome} \r \in N_0$ of $\ome
\backslash \{ 0\}$ into $\ome$ infinite
pieces.

Let us define by induction conditions $p^n
= \l p_0^n, 
\vtop{\offinterlineskip\hbox{$p_1^n$ } \hbox{$\scriptstyle\sim$}}
\r$ $(n < \ome)$.
Suppose that for every $n' < n$, $\l
p_0^{n'}, 
\vtop{\offinterlineskip\hbox{$p_1^{n'}$ } \hbox{$\scriptstyle\sim$}}
\r$ is defined.  
  Define
$\l p_0^n, 
\vtop{\offinterlineskip\hbox{$p_1^{n}$ } \hbox{$\scriptstyle\sim$}}
\r$. We consider $s_n$.
Let $s_n = \l k_{n_0} \nek k_{n\ell_n}\r$.
First define the coloring or the splitting
set of the condition.  Set $I_0 =  \{
\del_{k_{ni}} | i < \ell_n\}$.  Let us fill
the intervals between $\del_{k_{ni}}$'s.
We use  the colors of $J_{\l k_{n0}\r}$ to
color the ordinals below $\del_{k_{n0}}$
in a one-to-one fashion, i.e. no two
ordinals below $\del_{k_{n0}}$ are colored
by the same color.  Suppose that some fixed
function $f: J_{\l k_{n0}\r}
\leftrightarrow \del_{k_{n0}}$ was used all
the time for this purpose.

Then color $\del_{k_{n1}} \backslash
(\del_{k_{n0}} + 1)$ using $J_{\l k_{n0},
k_{n1}\r}$ as the set of colors and so on.
Let $\l I_i| i < \ome\r$ denote such
defined coloring of $\del_{k_{n\ell_n}} + 1$.

If for some $m < n$, $\l p_0^m, 
\vtop{\offinterlineskip\hbox{$p_1^{m}$ } \hbox{$\scriptstyle\sim$}}
\r$
decides that $\ualp$ and $s_n$ do not
contradict $p_0^m$, i.e. if $p_0^m = \l s,
A\r$, then $s_n$ is an end extension of
$s$ and $s_n \backslash \max s\subseteq A$.
Then we define $
\vtop{\offinterlineskip\hbox{$p_1^{n}$ } \hbox{$\scriptstyle\sim$}}
 =
\vtop{\offinterlineskip\hbox{$p_1^{m}$ } \hbox{$\scriptstyle\sim$}}
$ and
$p_0^n = \l s_n, A \backslash \max s_n\r$.

Suppose otherwise that there is no $m < n$
such that $\l p_0^m, 
\vtop{\offinterlineskip\hbox{$p_1^{m}$ } \hbox{$\scriptstyle\sim$}}
\r$ decides
$\ualp$ and $s_n$ do not contradict
$p_0^m$.

Consider $t_n = \l \l s_n, \ome \backslash
\max s_n \r$, $\l I_i | i < \ome\r\r$.  It
is clearly a condition in $\calP$  since
$I_0$ is finite and each $I_i$ contains at
most one element.  Also $t_n \in N_{\ell_n
+ 1}$.  If there is no $*$-extension of
$t_n$ deciding the value of $\ualp$ then
let $p^n = \l p_0^n, 
\vtop{\offinterlineskip\hbox{$p_1^{n}$ } \hbox{$\scriptstyle\sim$}}
\r = t_n$.
Suppose otherwise.  Then we choose inside
$N_{\ell_n + 1}$ a $*$-extension $t'_n$ of
$t_n$ deciding $\ualp$.  Notice that
$\ell th (t'_n)$ is countable ordinal and
that it is in $N_{\ell_n + 1}$.  So
$\ell th(t'_n) < \del_{\ell_n + 1}$.

We now want to extend $t'_n = \l t'_{n0},
\ut_{n1}\r$ to a condition $t^{''}_n$ of
the length $\del$.  The simplest approach
would be to use colors of $J_{s_n}$ and
color ordinals in $\del\backslash
\ell th(t'_n)$ in a one-to-one fashion.
But this may contradict 2.1.2 ((2b)$(*)$).
So we need to be more careful.  Consider 
$\{rng \ut_{n1} (\bet) | \bet < \ell th
t'_n\}$ as a collection of $\aleph_0$
almost disjoint subsets of $\ome$ by
2.1.2(2a).  Let $C \subseteq \ome$ be an
infinite set almost disjoint to all the
sets in this collection. 

Split $C$ into $\aleph_0$ disjoint infinite
sets $\l C_i | i < \ome\r$.  Let $\l c_{ij}
| j < \ome\r$ be an increasing enumeration
of $C_i$ $(i < \ome)$.  Let $\l \alp_i | 0
< i < \ome\r$ be an enumeration of
$\del\backslash \ell th t'_n$.  We now define
 a $\calP_U$-name of the color of
$\alp_i$ $(i < \ome)$ to be
$\{ \l s_n ^\smallfrown <r_1 \nek r_i \r,
\ome \backslash r_i \r, c_{ir_i} \r | \l
r_1 \nek r_i \r \in [\ome]^i, \ r_1 > \max
s_n\}$, i.e. the color of $\alp_i$ will be
$c_{ir_i}$, where $r_i$ is the $\ell_n + i$-th
element of the $\calP_U$-generic sequence.

Let $t^{''}_n$ be obtained from $t'_n$ by
adding these names to $\ut'_{n1}$.

\proclaim Claim 2.1.4.0.  $t^{''}_n$ is a
condition in $\calP$.

\pr We need to check 2.1.2 ((2b)$(*)$).  So,
let $I \subseteq \del$, $I \in V$ $p \ge
t^{''}_{n0}$ and a finite $J \subseteq
\ome$ be given.  First let us apply
(2b)$(*)$  to $\l p, 
\vtop{\offinterlineskip\hbox{$t'_{n}$ } \hbox{$\scriptstyle\sim$}}
\r , J$ and $I \cap \ell
th (
\vtop{\offinterlineskip\hbox{$t'_{n1}$ } \hbox{$\scriptstyle\sim$}}
)$.  There will be a finite
$a' \subseteq \ell th (
\vtop{\offinterlineskip\hbox{$t'_{n1}$ } \hbox{$\scriptstyle\sim$}}
)$ such
that for every finite $b \subseteq I
\backslash a'$ there is $p' \ge^* p$ such
that
$p' {\buildrul \calP_U \under \llvdash}
( \forall k \in J\  \forall \bet \in b\
(
\vtop{\offinterlineskip\hbox{$t_{n1}$ } \hbox{$\scriptstyle\sim$}}
 (\bet) \ne k))$ and $(\forall
\bet_1 \ne \bet_2 \in b\ 
\vtop{\offinterlineskip\hbox{$t_{n1}$ } \hbox{$\scriptstyle\sim$}}
 (\bet_1)
\ne 
\vtop{\offinterlineskip\hbox{$t_{n1}$ } \hbox{$\scriptstyle\sim$}}
 (\bet_2))$.

Now we like to extend $a'$ to some $a$ in
order to take care of the ordinals between
$\ell th (
\vtop{\offinterlineskip\hbox{$t'_{n1}$ } \hbox{$\scriptstyle\sim$}}
)$ and $\del$.  Let $p = \l
\vec s, A\r$ and $\vec s = s_n ^\smallfrown
\l r_1 \nek r_m\r$.  Then the colors of
each of the ordinals $\alp_1 \nek \alp_m$
are decided by $p$ and they are $c_{1r_1}
\nek c_{mr_m}$.  But it is important
that nothing is decided about the colors of the
$\alp$'s
in $(\del \backslash \ell th (
\vtop{\offinterlineskip\hbox{$t'_{n1}$ } \hbox{$\scriptstyle\sim$}}
))
\backslash \{ \alp_1 \nek \alp_n\}$. Set $a
= a' \cup \{ \alp_1 \nek \alp_m \}$.  Let
$b \subseteq I \backslash a$ be finite.
First pick  $p'\ge^* p$ witnessing (2b)$(*)$
for $b \cap \ell th (
\vtop{\offinterlineskip\hbox{$t'_{n1}$ } \hbox{$\scriptstyle\sim$}}
)$.  Let $p'
= \l \vec s, B\r$.  Consider the sets of
possible colors $
\vtop{\offinterlineskip\hbox{$t^{''}_{n1}$ } \hbox{$\scriptstyle\sim$}}
 (\alp)$ for
$\alp$'s in $b$.  These are finitely many
almost disjoint sets such that every
$
\vtop{\offinterlineskip\hbox{$t^{''}_{n1}$ } \hbox{$\scriptstyle\sim$}}
(\alp)$ for $\alp \in b
\backslash \ell th (
\vtop{\offinterlineskip\hbox{$t'_{n1}$ } \hbox{$\scriptstyle\sim$}}
)$ is
infinite.  Hence there is $k_0 < \ome$ such
that $
\vtop{\offinterlineskip\hbox{$t^{''}_{n1}$ } \hbox{$\scriptstyle\sim$}}
(\alp) \backslash k_0$ is
disjoint to every $
\vtop{\offinterlineskip\hbox{$t^{''}_{n1}$ } \hbox{$\scriptstyle\sim$}}
 (\bet)$
where $\alp \ne \bet \in b$, $\alp \in b
\backslash \ell th (
\vtop{\offinterlineskip\hbox{$t'_{n1}$ } \hbox{$\scriptstyle\sim$}}
)$.  Consider
$p^{''} = \l \vec s, B \backslash (1 + k_0
+ \max J) \r$.  We claim that $p^{''}$ is
as desired.  Let $\alp \in b \backslash
\ell th (
\vtop{\offinterlineskip\hbox{$t'_{n1}$ } \hbox{$\scriptstyle\sim$}}
)
$.  Find $i < \ome$ such
that $\alp = \alp_i$.  By the definition of
the set of possible colors for $\alp_i$, all
of them are above $\min (B \backslash (k_0
+ \max J))$, since $c_{ir} \ge r$ for every
$r < \ome$.  But this implies what we
need.\hb
$\bigsquare$ of the claim.

This complete the inductive definition of
$\l p^n | n < \ome\r$.

Let us now combine all $p^n$'s into
$\calP_U$-name $\uq1$.  Set $\uq1 = \{ \l
p_0^n, 
\vtop{\offinterlineskip\hbox{$p_1^n$ } \hbox{$\scriptstyle\sim$}}
\r | n < \ome\}$. 

\proclaim Claim 2.1.4.1.
$\l\l\ \r,
\ome \r, \uq1\r$  is a condition in
$\calP$.

\pr Let us check that $\l\l\ \r, \ome \r,
\uq1\r$
satisfies Definition
2.1.2.  The condition (1) of 2.1.2 is
trivial.

Let us show that $\l \l\ \r, \ome\r$ forces
 $\ell th (\uq1) = \del$, i.e. we
need to show that $\l\l\ \r, \ome \r$ forces 
every $\tau < \del$ to be colored but nothing
is colored above $\del$.  Let $\l s, A\r
\ge \l \l\ \r, \ome, \r$ and $\tau < \del$.
Pick $k < \ome$ so that $\del_k > \tau$.
Find $n < \ome$ such that $s_n$ is an end
extension of $s$ of  length $> k$ and $s_n
\backslash s \subseteq A$.  Then $\l s_n, A
\backslash \max s_n\r \ge \l s, A\r$ and by
the construction of $p^n$ it forces 
$\tau$ to be colored.  Suppose now that $\tau
\ge \del$.  In $V^{\calP_U}, \uq1$ is
interpreted as a coloring of ordinals less
than $\del$ only, by its construction.  So
we are completely free about the color of
$\tau$.

Let us now check  the condition 2.1.2
((2b)$(*)$). So, suppose that $I \subseteq
\del$, $I \in V$, $p'_0 = \l \vec s, A\r
\ge \l\l\ \r, \ome\r$ and a finite $J
\subseteq \ome$ are given.  If for some $B
\subseteq \ome$ $\l \l \vec s, B \r,
\uq1\r$ decides $\ualp$ then we are
exactly in the situation of Claim 2.1.4.0.
So suppose otherwise.  But then nothing
nontrivial was done at stage $n$ with $s_n
= \vec s$.  So the set $a = \{
\del_{k_{ni}} | i < \ell_n\}$ will work,
where $s_n = \l k_{n0} \nek
k_{n\ell_n}\r$.\hb
$\bigsquare$ of the claim.

So $\l\l\l \ \r, \ome\r, \uq1\r \in \calP$.
Let us show that it is as desired.  We need
to define a $\calP_U$-name $\ubet$ such
that $\l\l\l\ \r, \ome\r, \uq1\r
\ldash\ualp = \ubet$.  Set $\ubet = \{\l 
p_0^n, \check \gam_n \r | n < \ome, \l
p_0^n, 
\vtop{\offinterlineskip\hbox{$p_1^n$ } \hbox{$\scriptstyle\sim$}}
\r \ldash \ualp = \check
\gam_n\}$.

Let $\l r_0, \ur1\r$ be an extension of
$\l\l\l\ \r , \ome\r, \uq1\r$ forcing
$''\ualp = \check \gam ''$ for some ordinal
$\gam$.  Consider $r_0$.  Let $r_0 = \l s,
A\r$.  Then for some $n < \ome$, $s = s_n$.
Hence $r_0$ is compatible with $p_0^n$.  Then
by
the construction, a $*$-extension
$t'_n$ of $t_n$ deciding the value of
$\ualp$ to be some $\gam_n$ was picked
and $r_0 \ldash^{''} \ur1$ is an end
extension of $t_{n1}^{'\ ''}$.  Hence $\gam_n =
\gam$.  This completes the proof of the
lemma.\hb
$\bigsquare$

\proclaim Lemma 2.1.5.  Let $\uf$ be a
$\calP$-name of a function from $\ome$ into
ordinals.  Then there exist 
$\calP_U$-names $\ug$ and $q_1$ such that
$$\l\l\ \r, \ome\r, \uq1\r \ldash \ug =
\uf\ .$$

\pr Let $\l N_n | n < \ome\r$ be an
increasing sequence of countable elementary
submodels of a large enough portion of the
universe such that

(a) $\calP, \uf \in N_0$

(b) $N_n \in N_{n+1}$ for every $n < \ome$.

Let $N = \bigcup\limits_{n < \ome} N_n$,
$\del_n = N_n \cap \ome_1$ and $\del =
\bigcup\limits_{n <  \ome} \del_n = N \cap
\ome_1$.  Let $\l s_n | n < \ome\r$ and $\l
J_{\vec n} | \vec n \in [\ome]^{< \ome} \r$
be as in Lemma 2.1.4.

We define by induction a sequence of
conditions $p^n = \l p_0^n, 
\vtop{\offinterlineskip\hbox{$p_1^n$ } \hbox{$\scriptstyle\sim$}}
\r\ (n <
\ome)$.  Suppose that for every $n'< n$ $\l
p_0^{n'}, 
\vtop{\offinterlineskip\hbox{$p_1^{n'}$ } \hbox{$\scriptstyle\sim$}}
\r$ is defined. Define $\l p_0^n,
p_1^n\r$.  We consider $s_n$.
Let $s_n = \l k_{n0} \nek k_{n \ell_n} \r$.
Find $m < n$ such that $s_m = \l k_{n0}
\nek k_{n, \ell_n-1}\r$.  Assume as an
inductive assumption that $p^m \in
N_{k_n\ell_{n-1} + 1}, \ell th (p^m) \ge
\del_{k_{n \ell_n - 1}}$ and for every $i <
\ell_n$, $p^m$ decides the value of
$\uf(i)$.

Pick an infinite set $C_{s_m} \subseteq
\ome$ which is almost disjoint with each
set of possible colors $
\vtop{\offinterlineskip\hbox{$p_1^m$ } \hbox{$\scriptstyle\sim$}}
 (\bet)$ for
$\bet< \ell th (p^m)$.  Split $C_{s_m}$
into $\ome$ disjoint infinite sets $\l
C_{s_m, \vec n} | \vec n \in [\ome]^{<
\ome}$, $\vec n$ is an end extension of
$s_m \r$.  At stage $n$ we are going to use
$C_{s_m, s_n}$.  Split this set again into
$\aleph_0$ disjoint infinite sets $\l C_i |
i < \ome\r$.  Let $\l c_{ij} | j < \ome \r$
be an increasing enumeration of $C_i $ $(i
< \ome)$.  Let $\l \alp_i | 0 < i < \ome\r$
be an enumeration of $\del_{k_{n \ell_n}}
\backslash \ell th (p^m)$.  Define a
$\calP_U$-name of the color 
of $\alp_i$ $(i < \ome)$ to be $\{ \l s_n
^\smallfrown < r_1 \nek r_i \r, \ome
\backslash r_i \r, c_{i r_i} | \l r_1 \nek
r_i\r \in [\ome]^i, r_1 > \max s_n\}$.  So
$\{ c_{ij} | j < \ome\}$ will be the set of
possible colors of $\alp_i$ and if $r_i$ is
$\ell_n + i$-th element of the
$\calP_U$-generic sequence then the color
of $\alp_i$ will be $c_{ir_i}$.  Color 
$\del_{k_n\ell_n}$ into $0$ color.  Let
${\buildrul \sim \under p}^{n'}$ be obtained from ${\buildrul \sim \under p}^m$ by
doing this additional coloring.  Assuming
$C$ and all the portions and enumerations
used are in $N_{k_{n \ell_n} + 1}$, we obtain
that ${\buildrul \sim \under p}^{n'} \in N_{k_{n \ell_n} + 1}$.
Repeating Claim 2.1.4.1, one can show that
$\l s_n, \ome \backslash \max s_n,
{\buildrul \sim \under p}^{n'}\r$ is a condition in $\calP$.  Now
apply Lemma 2.1.4 inside $N_{k_{n\ell_n} +
1}$ to this condition in order to decide
$\uf (\ell_n)$.  Let $p_0^n = \l s_n,
\ome-\max s_n\r$ and $
\vtop{\offinterlineskip\hbox{$p_1^{n}$ } \hbox{$\scriptstyle\sim$}}
$ be such
a condition as produced by Lemma 2.1.4.

This completes the inductive definition.

 Now define $\uq1 = \{ \l p_0^n, 
\vtop{\offinterlineskip\hbox{$p_1^{n}$ } \hbox{$\scriptstyle\sim$}}
 \r
| n < \ome\}$.  It is enough to show that
$\l\l\l\ \r, \ome \r, \uq1\r\in \calP$.

As in Claim 2.1.4.1, $\ell th (\uq1) =
\del$.  So let us check (2b)$(*)$ of 2.1.2.
Suppose $I \subseteq \del$, $I \in V$,
$p'_0 \ge \l\l\  \r, \ome \r$ and $J \subseteq
\ome$ are given.  Without loss of
generality $I = \del$.  Let $p'_0 = \l s,
A\r$ and $s = \l k_0 \nek k_\ell\r$.  For
some $n < \ome$, $s = s_n$.  Consider first
$\l\l k_0\nek k_\ell\r, \ome\backslash
 k_\ell\r, 
\vtop{\offinterlineskip\hbox{$p_1^{n}$ } \hbox{$\scriptstyle\sim$}}
\r$.  It is a condition in
$\calP$ of the length $< \del_{k_\ell +
1}$.  Let $a' \subseteq \ell th (
\vtop{\offinterlineskip\hbox{$p_1^{n}$ } \hbox{$\scriptstyle\sim$}}
)$ be
the set witnessing (2b)$(*)$ for this
condition.  Set $a = a' \cup \{ \del_{k_i} |
i \le \ell\}$.  Suppose that $b \subseteq
\del \backslash a$ is finite.  Consider the
set $C_{s_n}$ which was used in the
definition of $p_n$ (or $\uq1$) and which
is almost disjoint to each set of possible
colors $
\vtop{\offinterlineskip\hbox{$p_1^{n}$ } \hbox{$\scriptstyle\sim$}}
 (\bet)$ for $\bet < \ell th
(p^n)$.  It was split into $\l C_{s_n, \vec
n} | \vec n \in [\ome]^{< \ome}, \vec n$ is
an end extension of $s_n\r$.

Find $m_0$, $\max s_n < m_0 < \ome$ such that for
every finite end extension $\vec n$ of $s_n$
with $m_0 \le \max \vec n$ $C_{s_n \vec n}$
is disjoint to $J$ and to $
\vtop{\offinterlineskip\hbox{$p_1^{n}$ } \hbox{$\scriptstyle\sim$}}
(\bet)$ for $\bet \in b \cap \ell th
(
\vtop{\offinterlineskip\hbox{$p_1^{n}$ } \hbox{$\scriptstyle\sim$}}
)$.  Let $m_1 > m_0$ be such that
$\del_{m_1} > \max b$.  Let $\l s, A'\r$ be
witnessing (2b)$(*)$ for $b \cap \ell th
(
\vtop{\offinterlineskip\hbox{$p_1^{n}$ } \hbox{$\scriptstyle\sim$}}
)$.  Set $\tilA = A' \backslash
m_1$ and $p^{''}_0 = \l s, \tilA\r$.  Then
$ p_0^{''} {\buildrul  \calP_U\under \llvdash}
(\forall k \in J \ \forall \bet \in b\  (\uq1
(\bet) \ne k))$ and $(\forall \bet_1 \ne
\bet_2 \in b\  \uq1 (\bet_1) \ne
\uq1(\bet_2))$.

The first part is clear since the possible
colors of the ordinals above $\ell th
(p^n)$ and below $\max b + 1$ are disjoint
to $J$.  For the second notice that $\max
b + 1$ is forced to be below the next
$\del_k$ where the value of $\uf$ is
decided.  So ordinals in the interval
$[\ell th (p^n), \max b + 1]$ are colored
using  a disjoint partition of $C_{s_n \vec n}$
for some $\vec n$.  So the colors of the elements of
$b \backslash \ell th (p^n)$ are different
(or are forced to be such).  By the choice
of $m_0$, these colors are disjoint to the
sets of possible colors of $
\vtop{\offinterlineskip\hbox{$p_1^{n}$ } \hbox{$\scriptstyle\sim$}}
 (\bet)$
for $\bet \in b \cap \ell th (p^n)$. So $
p^{''}_0$ will be as desired.

This completes the checking of (2b)$(*)$ and
hence the proof of the lemma.\hb
$\bigsquare$

The following is now immediate.

\proclaim Corollary 2.1.6.  The forcing $\l
\calP, \le \r$ preserves the cofinalities.

\proclaim Corollary 2.1.7.  $V^\calP
\models $GCH.

Let us now turn to the proof of 2.1.

\medskip

\noindent
{\bf Proof of Theorem 2.1:} Force with $\l
\calP, \le \r$.  Let $G \subseteq \calP$ be
a generic set.  Define $V_1 = V[G]$.  Now force
 a \C real $r$ over $V_1$ 
(let us deal with \C reals; the random real
case is similar).  We want to produce
$\aleph_1$-\C reals over $V$. Fix a
sequence $\l f_i | i < \ome_1\r \in V$ of
almost different functions from $\ome$ to
$\ome$, i.e. for  every $i \ne j$ $|\{ n |
f_i(n) = f_j(n) \}| < \aleph_0$.

Split $r$ into $\ome$-\C reals over $V_1$.
Denote them by $\l r_{n,m} | n,m < \ome
\r$.  Now define  the $\ome_1$-reals $\l s_{n,i}
| n < \ome, i < \ome_1\r$ as follows:
$$s_{n,i} (k) = r_{n,f_i (k)} (0)$$
where $k < \ome$.

Let $\l I_n | n < \ome \r$ be the portion
of $\ome_1$ produced by $G$.  We spread $\l
s_{n,i} | i < \ome_1\r$ inside $I_n$.  More
precisely, for every $\alp < \ome_1$ let
$t_\alp = s_{n, i}$ if $\alp \in I_n$ and
$\alp$ is the $i$-th element of $I_n$.

\proclaim Claim 2.1.8. $\l t_\alp | \alp <
\ome_1 \r$ are $\aleph_1$-\C generic reals
over $V$.

\pr Suppose otherwise.  Let $\l\l p_0,
\upp1\r, q\r \in \calP \times$ \C
$(\aleph_0 \times \aleph_0, 2)$ and $D
\subseteq$ \C $(\aleph_1, 2)$ dense open,
 $D \in V$ be witnessing this. 

Let $D^* \subseteq D$, $D^* \in V$ be a
maximal antichain generating $D$.  Set $I =
\{ \alp < \ome_1 | \exists t \in D^*\ \alp \in
\dom t\}$. Let $\alp = \sup I$.  Clearly
$\alp < \ome_1$.  Let us assume that $\ell
th (\upp1) \ge \alp$ otherwise just extend
it to such a condition.  Now let us apply
(2b)(*) of 2.1.2 to $\l p_0, \upp1\r$, $I$
and $J = \{ n | \exists m \ (n,m) \in \dom
q\}$.  Let $a$ be a finite subset of $I$
given by (2b)(*).  For every $\nu \in a$
find unique $n (\nu)$ and $i(\nu)$ such
that $\nu \in I_{n(\nu)}$ and it is
the $i(\nu)$-th element of $I_{n(\nu)}$.  If
$n(\nu) \not \in \dom (\dom q)$ or $n(\nu)
\in \dom (\dom q)$ but for no $k < \ome$
$(n(\nu), f_{i(\nu)} (k) ) \not\in \dom
(q)$, then $q$ carries no information on
$t_\nu = s_{n(\nu), i(\nu)}$.  Suppose
otherwise.  Then for every $k < \ome$ with
$(n(\nu), f_{i(\nu)} (k)) \in \dom q$ the
value of $t_\nu (k)$ is decided  by
$q(n(\nu), f_{i(\nu)} (k))$.  Let $s= \{ \l
\nu, k, j\r | \nu \in a, k<\ome, j < 2$ and
$(n(\nu), f_{i(\nu)} (k), j) \in q\}$.  It
is actually all the information provided by
$q$ about $t_\nu$'s with $\nu \in a$.  Then
$s \in {\rm Cohen}\ (\aleph_1)$.

There is $k_0 < \ome$ such that for every
$k > k_0$ for every $\nu_1, \nu_2 \in a$,
$f_{i(\nu_1)} (k) \ne f_{i(\nu_2)} (k)$.
For every $\nu \in a$ and every $k \le k_0$
if $\l \nu, k, 1\r \not\in s$ then add $\l
\nu, k, 0\r$ to $s$.  Let $s^*$ be a
condition obtained from $s$ this way.

The reason for defining $s^*$ is to avoid
possible collisions.  Remember that
$s_{n,i} (k)$ was defined to be $r_{n,
f_i(k)} (0)$.  Hence, if $f_{i_2} (k) =
f_{i_2} (k)$ then it is necessary that $s_{n, i_1} (k)
= s_{n, i_2} (k)$.

Now choose some $s^{**} \ge s^*$ such that
$s^{**} \in D$ and $\dom (\dom s^{**})
\subseteq I$.  Set $b = \dom(\dom s^{**})
\backslash a$.  Using (2b) (*) of 2.1.2,
find $p'_0 \ge^* p_0$ forcing
$$\eqalignno{
&''\big(\forall k \in J\ \forall \bet \in b\
\ (\upp1 (\bet) \ne k) \big)''\cr
\noalign{\hbox{and}}
&''\big(\forall \bet_1 \ne \bet_2 \in b\ \
(\upp1(\bet_1) \ne \upp1 (\bet_2)\big)''\
.\cr}$$
We extend $p'_0$ further to a condition
$p^{''}_0$ deciding all the colors of
elements of $a \cup b$.  Now there is no
problem to extend $q$ to some $q^*$ such
that above $\l p^{''}_0,
\upp1\r$ in $V^\calP$ $q^* \ldash$'' for every $\nu
\in \dom(\dom s^{**}), s^{**} (\nu)$ is a
subfunction of $t_\nu$".

So $\l\l p^{''}_0, \upp1\r, q^*\r \ge \l \l
p_0, \upp1 \r, q\r$ forces  $D$  not to be
a counterexample.  Contradiction.\hb
$\bigsquare$

A similar construction may be used to give
another proof of a result of
B. Velickovic and H. Woodin.

\proclaim Theorem {\rm [Ve-Wo]}.  There is
a pair $(W, W_1)$ of generic extensions of
$L$ such that $W \subseteq W_1$  $\aleph_1^W
= \aleph_1^{W_1}$ and $W_1$ contains a
perfect set of $W$-reals which is not in
$W$.

Namely we get a slightly stronger result:

\proclaim Theorem 2.2.  There is a pair $(W,
W_1)$ of generic cofinality preserving
extensions of $L$ such that $W \subseteq
W_1$ and $W_1$ contains a perfect set of
$W$-reals which is not in $W$.

\pr We force over $V^\calP$ a \C real and a
perfect set of \C reals (or just a perfect
set of \C reals).

Now we spread this perfect set on the set
of ordinals colored by 0 and using the \C
real as in Theorem 2.1 produce
$\aleph_1$-reals indexed by the rest of
colors.  Denote this combined sequence by
$\l t_\alp  |\alp < \aleph_1\r$.  Repeating
the argument of 2.1, one can see that $\l
t_\alp | \alp < \ome_1\r$ are forming
$\aleph_1$-Cohen generic reals over $V$.
Set $W = V [t_\alp | \alp < \ome_1]$, $W_1
= V^{\calP * \ \hbox{perfect set of Cohen
reals}}$.\hb
$\bigsquare$

\vskip 1truecm

\references{50}

\ref{[Fn]} D. Fremlin, Real-valued
measurable cardinals, in Set Theory of the
Reals, H. Judah, ed., Israel Math. Conf.
Proceedings (1993), 151-305.

\ref{[Gi-Ma]} M. Gitik and M. Magidor, The
SCH revisited, in Judah-Just-Woodin, eds,
Set Theory of the Continuum, MSRI
publications \# 26, Springer 1992, 243-279.

\ref{[Gi-Ma]} M. Gitik and M. Magidor,
Extender based forcings, JSL, 59 (1994),
445-460.

\ref{[Ma]} M. Magidor, Changing cofinality
of cardinals, Fund. Math. XCIX (1978),
61-71.

\ref{[Ve-Wo]} B. Velickovic and H. Woodin,
Complexity of the reals of inner models of
set theory.

\end
\vtop{\offinterlineskip\hbox{$I_{p,n}$ } \hbox{$\scriptstyle\sim$}}